\documentstyle[12pt]{article}
\begin{document}

\title{Linear Connections on Graphs } 
\author{ Sunggoo Cho\\Department of Physics, Semyung University, \\
 Chechon, Chungbuk 390 - 230, Korea \\ Kwang Sung Park \\ Department
 of Mathematics, Keimyung University, \\
Taegu 705 - 701, Korea \\}
\maketitle

\begin{abstract}
In recent years,  discrete spaces
such as  graphs attract much attention as models 
for physical spacetime or as models for testing the spirit of 
non-commutative geometry. 
In this work,  we 
 construct the  differential algebras for graphs by extending
the work of Dimakis et al and
discuss linear connections and curvatures on  graphs.
Especially, we calculate connections and curvatures explicitly for the 
general nonzero torsion case.  There is a metric, but no metric-compatible
connection in general except the complete symmetric graph with two vertices.   
 \\

\vspace{2cm}
\vspace{2cm}

PACS NO: 02.40.+m, 04.60.+n
\vspace{0.5cm}

SMHEP 96-12

\end{abstract}

\newpage

{\bf \Large I. INTRODUCTION }
\vspace{0.5cm}

In the last few years, there has been a rapid increase of interest
in  non-commutative geometry.
Non-commutative geometry is the geometry of quantum spaces, which
are generalized spaces replacing classical smooth
manifolds$^{1}$.
 A quantum space is an associative algebra,  usually non-commutative,
with possibly more structures on it$^{1,2}$.
Even though non-commutative geometry is developed usually for
non-commutative algebras, however, there can be interesting commutative
 quantum spaces such as discrete sets, which we are mainly concerned 
with in this work.
In fact, discrete spaces attract much attention in recent years as models 
for physical spacetime or as models for testing the spirit of 
quantum spaces$^{3-6}$.  
Differential calculi on commutative algebras 
 have already 
been investigated by some authors$^{6-8}$.

One of the motivations for studying quantum spaces 
is to  understand the small scale structure
of spacetime and quantum gravity.
 There are some efforts to try to understand  
gravity in the framework of the
non-commutative
version of Riemannian geometry (see, e.g., Ref.\ 9, 10 ).
 However, non-commutative Riemannian geometry is not a 
straightforward extension of the ordinary one and  
 several definitions of a connection have been considered.
In particular,   definitions of a linear connection 
which make use of the bimodule 
structure of differential forms have been 
suggested for  projective modules$^{11}$ and
for differential calculi based on derivations$^{12,13}$.
For discrete sets, one can define a linear connection
 using a group structure supplied to the discrete sets$^{14}$. 
 Recently, a more general definition has been proposed by Mourad$^{15}$
of a linear connection, 
which is  also our main concern in this work.

In this work, we shall extend the formulation of 
Dimakis et al$^{6}$ on graphs in order to make them suitable for
calculation of linear connections proposed by Mourad and curvatures.  
 In sec.\ II, we shall have a review of 
  the universal differential algebra
 necessary in the sequel.
And  differential caculi on a graph$^{6}$ is extended toward
 the construction of the differential algebra for the graph.
In sec.\ III, we shall calculate explicitly Mourad's linear connections 
and curvatures on graphs.
 Nonzero torsion connections, bilinear curvatures and metrics shall be
discussed.
\vspace{1cm}

{\bf \Large II. DIFFERENTIAL ALGEBRAS FOR GRAPHS }
\vspace{0.5cm}

{\bf \large A. Universal differential algebra}
\vspace{0.5cm}

Let $A$ be an associative algebra with 1 over the field {\bf C} of complex 
numbers.
Let a direct sum of vector spaces $\Omega = \oplus _{n=0 }^{\infty}
\Omega ^{n}$ 
 be a  differential complex so that there are 
homomorphisms $d$'s 
\begin{equation}
  \cdots \rightarrow \Omega ^{n-1} \stackrel{d}{\rightarrow }       
\Omega ^{n} \stackrel{d}{\rightarrow} \Omega ^{n+1} \rightarrow \cdots
\label{eq:1}      
\end{equation}
such that $d^{2} = 0$.  The homomorphism $d$ is usually called the
 differential operator of
the complex $\Omega $. 
 If there is a gradation-respecting multiplication
 $\cdot $ in $\Omega $ so that $\Omega $ is an algebra over {\bf C} 
 and the homomorphisms $d$'s satisfy the Leibniz rule
\begin{equation}
 d(\omega \cdot \omega ') = (d\omega)\cdot \omega ' +
 (-1)^{n}\omega \cdot d\omega '  
\label{eq:2}      
\end{equation}
where $\omega \in \Omega^{n} $, then
the algebra is called a differential algebra  over
 {\bf C}. 

There is an important example for a differential algebra, which is 
 constructed from an associative algebra $A$  with 1 over {\bf C}
as follows (See, e.g., Ref.\ 16, 17).
Let $\Gamma  $ be
 an $A$-bimodule.
Let   $d : A \rightarrow \Gamma $ be a linear map such that 
for any $a,b \in A $, 
  $d(ab) = (da)b +a\,db $.  
If every element of $\Gamma $ is of 
the form $\sum_{k}a_{k}db_{k} $ where $a_{k}, b_{k} \in A $, then 
$(\Gamma,  d)$ is said to be
 a 1st order differential calculus over $ A$. 
We have a special 1st order differential calculus over $A$.
Let $m : A \otimes _{{\bf C}} A
 \rightarrow A $ be the multiplication map in $A$ such that $m(a \otimes b)=
ab$.  Let $\Omega^{1} \equiv ker \, m $ and let $d : A \rightarrow \Omega^{1}$
be a map defined by $da \equiv 1\otimes a - a \otimes 1 $.  Then
$\Omega^{1} $ is an $A$-bimodule and $(\Omega^{1}, d )$ is
 a 1st order differential
calculus over $A$. 
 Since 
 $A \otimes _{{\bf C}} A$ also 
carries an $A$-bimodule structure and
 the multiplication map $m : A \otimes _{{\bf C}} A \rightarrow A $ is a
 bimodule homomorphism, we have an exact sequence of bimodule
 homomorphisms
\begin{equation}
 0 \rightarrow \Omega^{1} \stackrel{\imath }{\rightarrow} A 
\otimes_{{\bf C}}A
\stackrel{m}{\rightarrow} A \rightarrow 0 . 
\label{eq:3}      
\end{equation}
It is well known that the
 1st order differential algebra  $(\Omega^{1}, d ) $ over $A$
 is universal, i.e. if there is another 1st order differential calculus
$(\Gamma, \delta )$ over $A$, then there exists a unique $A$-bimodule
 homomorphism $\phi : \Omega^{1} \rightarrow \Gamma $ such that
\begin{equation}
\delta = \phi \circ d.
\label{eq:4}
\end{equation}
Every 1st order differential calculus 
$(\Gamma ,  \delta )$ over $A$ is isomorphic, as $A$-bimodules,
 to a quotient of $\Omega^{1} $
by the  $A$-submodule $ker \, \phi $.
 
We can extend the 1st order differential calculus $\Omega^{1} $ 
to higher orders:
 Let $\Omega ^{n} \equiv  \{ \rho \in A \otimes _{{\bf C}}
\cdots \otimes _{{\bf C}} A \equiv A^{\otimes _{{\bf C}}(n +1)}
\mid 
 m_{i} \rho = 0 $ for all $i = 1,\cdots , n \}$ where $m_{i}$
is the multiplication acting in the $i, (i+1)$th place.
Then  $\Omega^{n} $ is an $A$-bimodule and
  $\Omega ^{n} = {\em span}\{a_{0}da_{1} \otimes \cdots \otimes da_{n} 
                   \mid a_{i} \in A \}$\footnote{ $\! ^)$We shall
 use the same notation
$ \otimes  $ for the two kinds of tensor products 
$ \otimes _{{\bf C}}  $, $ \otimes _{A} $  if
there is no confusion. }$\! ^)$
   $ = \Omega^{1}\otimes_{A} \cdots \otimes_{A} \Omega^{1} \equiv
 (\Omega^{1})^{\otimes _{A}n}. $
From these $\Omega^{n} $'s, we can construct a differential algebra.
Let $\Omega (A) \equiv \oplus_{n=0}^{\infty }\Omega^{n} $ 
where $\Omega^{0} \equiv A $. And we define a multiplication $\cdot $
 in $\Omega(A) $ as follows:

$(i)$ Let $a_{0}da_{1}\otimes da_{2} \cdots \otimes da_{n} \equiv 
      (a_{0}, a_{1}, a_{2}, \cdots, a_{n}) $  
and define a multiplication of two elements in $\Omega^{n}$ and $\Omega^{m-1}$
to get an element in $\Omega^{n+m-1} $by 
\begin{eqnarray}
&&(a_{0}, a_{1},  \cdots,  a_{n}) \cdot ( a_{n+1}, \cdots , a_{n+m}) 
                                                              \nonumber \\
& & = \sum_{i=0}^{n}(-1)^{n-i}
  (a_{0}, a_{1}, \cdots ,a_{i}a_{i+1}, \cdots ,
                 a_{n}, a_{n+1}, \cdots , a_{n+m}). 
\label{eq:5}  
\end{eqnarray}

$(ii)$ Extend this multiplication to the whole of $\Omega(A) $ using the
 distributive rule in an obvious manner.

$(iii)$ Extend $d$ to the whole space $\Omega (A)$ as follows: \\
\begin{equation}
d(a_{0}, a_{1}, \cdots , a_{n})
  \equiv (1, a_{0}, a_{1}, \cdots , a_{n})
\label{eq:6}      
\end{equation}
and 
\begin{equation}
d(1,a_{0},a_{1}, \cdots , a_{n}) \equiv 0. 
\label{eq:7}      
\end{equation}
Then $\Omega (A) $ is a differential algebra. This differential
 algebra is universal.
\\    
It is straightforward to see that if $a_{i}a_{i+1} = 0 $ 
for $i = 0,1, \cdots, n-1 $, then
\begin{equation}
  (a_{0}, a_{1}, \cdots , a_{n}) = a_{0} \otimes a_{1} \otimes
   \cdots \otimes a_{n},
\label{eq:8}         
\end{equation}
\begin{equation}
(a_{0}, a_{1},  \cdots,  a_{n}) \cdot ( a_{n+1}, \cdots , a_{n+m}) 
  =  (a_{0}, \cdots , a_{n-1},a_{n} a_{n+1}, a_{n+2}, \cdots , a_{n+m}) 
\label{eq:9}        
\end{equation}
and
\begin{equation}
 d(a_{0}, a_{1}, \cdots, a_{n})
= \sum_{q=0}^{n+1}(-1)^{q}a_{0}\otimes
   a_{1} \otimes \cdots \otimes a_{q-1} \otimes 1 \otimes a_{q}\otimes \cdots
 \otimes
 a_{n} . 
\label{eq:10}      
\end{equation} 
\vspace{0.5cm}

{\bf \large B. Differential algebras on graphs}
\vspace{0.5cm}

Let $\cal V$ be a set of $N$ points $x_{1}, \cdots , x_{N}\,\,
 (N < \infty )$.
As in Ref.\ 6, let $A$ be the algebra of complex functions on ${\cal V}$ with
 $(fg)(x_{i}) = f(x_{i})g(x_{i}) $.  Let $e_{i} \in A $ be defined by

\begin{equation} 
e_{i}(x_{j})=\delta _{ij}.   \label{eq:11}
\end{equation}
Then it follows that
\begin{equation}
e_{i}e_{j} = \delta_{ij}e_{i}, \hspace{1cm} \sum_{i}e_{i} = 1    \label{eq:12}
\end{equation}
and each $f \in A$ can be written as $f = \sum_{i}f(i)e_{i}$
where $f(i) = f(x_{i}) \in {\bf C} $.
It is obvious that $A$ is not only a commutative algebra with 1
but also
an $N$-dimensional complex vector space.

Now let us introduce the universal differential algebra $\Omega (A)$ and
the differential operator $d$ as in the previous section.
Then the differentials satisfy the following relations.
\begin{equation}
e_{i}de_{j} = -(de_{i})e_{j} + \delta_{ij}de_{i}  \label{eq:13}
\end{equation}
and
\begin{equation}
\sum_{i}de_{i} = 0,\hspace{1cm} d1 = 0.  \label{eq:14}
\end{equation}
The universal 1st order differential calculus $\Omega^1$ is generated by
 ${\cal B}\equiv
 \{ e_{i}de_{j} \mid i, j = 1, 2, \cdots , \,\,\, N ( i \ne j )\} $ 
as an $A$-bimodule.
In this work, we note that $\Omega^1$ is a finite-dimensional
 complex vector space 
with having the generators  
$\cal B $
as a basis.

Similarly,
for $n \geq 2$,  $\Omega^n$ is not only an  $A$-bimodule generated by
$ \{ e_{i_{1}}de_{i_{2}} \otimes \cdots
 \otimes de_{i_{n+1}} \mid i_{k} = 1, 2, \cdots ,
 N \,\,\,(i_{k} \ne i_{k+1}) \}   $
but also  a complex vector space with having the generators as a basis.
In the universal differential algebra $\Omega(A)$ of $ A$,
 the multiplication $\cdot $ in Eq.(\ref{eq:9}) yields

\begin{eqnarray}
&&(e_{i_{1}}de_{i_{2}} \otimes  \cdots \otimes de_{i_{r}})   \cdot 
 ( e_{j_{1}}de_{j_{2}} \otimes \cdots
 \otimes de_{j_{s}})  \nonumber \\
&& =   
 e_{i_{1}}de_{i_{2}} \otimes \cdots \otimes d(e_{i_{r}}e_{j_{1}}) \otimes \cdots
    \otimes de_{j_{s}} \label{eq:15}         \\
 & &=  
  \delta_{i_{r}j_{1}}e_{i_{1}}de_{i_{2}} \otimes \cdots  \otimes
de_{i_{r-1}} \otimes de_{j_{1}} \otimes de_{j_{2}} \otimes \cdots 
\otimes de_{j_{s}} .
 \nonumber
\end{eqnarray}
Thus, effectively,
the multiplication $\cdot $ is   the same as the tensor product
$\otimes _{A} $, which is crucial in the sequel together with the fact that
the differential operator $d$ satisfies Eq.\ (\ref{eq:10}).

A graph is  a set of vertices which are 
interconnected by a set of edges$^{18}$. 
 Graphs are assumed to be connected,
i.e. two arbitrary vertices can be connected by a sequence of 
consecutive edges.
A complete graph is a graph 
   for which every pair of distinct vertices is connected by one edge.
In some applications, it is natural to assign a direction to 
each edge of a graph. 
 Diagrammatically, the direction of each edge 
is represented by an arrow. 
A directed graph (or digraph ) is a graph augmented in this way.
A symmetric digraph is a digraph in which any connected pair 
of distinct vertices  is connected in both directions.
In this work, we shall be concerned only with digraphs and hence every
graph is a digraph unless otherwise stated.

Now let us regard the points $x_{i}$ of $\cal V $
as vertices and a generator $e_{i}de_{j} (i \ne j )$ in ${\cal B}$
 as an arrow from $x_{i}$ to 
$x_{j}$.  Then we can associate a graph with a set $\cal V$ of points 
and a subset ${\cal S}$ of  ${\cal B}$. 
 From now on, we shall denote
a graph by $(A, K^1 )$, where $K^1 $ is an $A$-bimodule generated 
by the subset ${\cal S}$.
 Accordingly,
 $(A, \Omega^{1})$ 
represents a complete symmetric graph 
since the generators in ${\cal B}$ correspond to all arrows connecting
any two vertices.
 In this case, the graph
$(A, K^1 )$ is said to be a subgraph of $(A, \Omega^1 )$.
A subgraph $(A, K^1)$ is obtained by deleting some of the arrows in a 
complete symmetric graph $(A, \Omega^1)$.
It is obvious that $\Omega^1 $ is the direct sum of $K^1$ and 
its complement which is an $A$-bimodule generated by the set 
${\cal B} - {\cal S}$.
We define an $A$-bimodule homomorphism $\phi_{1} : \Omega^1 \longrightarrow
K^1 $ to be the projection map.  If we define a map
$\delta : A \longrightarrow K^1 $ by $\delta = \phi_{1} \circ d $, then
it is straightforward to see that $( K^1, \delta )$ is a 1st order 
differential calculus over $A$ 
with an observation
$e_{i}de_{j} = \phi_{1}(e_{i}de_{j})
             = e_{i}\phi_{1}(de_{j})
             = e_{i}\delta e_{j}$ 
for $e_{i}de_{j} \in K^1 $.

Now we shall construct a differential algebra from $K^{1}$.
We  
define $K^{n}$ to be the quotient space:
\begin{equation}
 K^n  \equiv \frac{\Omega^{n}}{ <d(ker \phi_{n-1})>}  
\label{eq:015}
\end{equation}
     for $ n = 2, 3, 4, \cdots $,
where $\phi_{n-1} $ is the projection map
  from $\Omega^{n-1} $ to $K^{n-1}$ and the bracket $ < \,X \, > $
 means the $A$-bimodule
generated by the set $X$. So
$ <d(ker \phi_{n-1})> $ is the 
      $A$-bimodule generated
 by the subspace $ d(ker \phi_{n-1}) $ of $\Omega^{n} $, and
 $ker \phi _{n} = <d(ker \phi_{n-1})> $. 
We
define $ \delta : K^{n} \longrightarrow K^{n+1} $ by
$ x \,+ < d(ker \phi_{n-1})> \mapsto \phi_{n+1}(dx) $
for $ x \in \Omega^{n} $ and for $n = 1, 2, \cdots $.
  This is well-defined since
if $ x \,+ <d(ker \phi_{n-1})>\, = \,y\, + <d(ker \phi_{n-1})> $,  then
$d( x - y ) \in \,\, <d(ker \phi_{n})> $ and hence  
$\delta (x \,+ <d(ker \phi_{n-1})> ) = \delta (y \,+ <d(ker \phi_{n-1})>)$.

Now let us investigate the form of   generators 
for $ker \phi_{n} $ $(n \geq 2)$
since it is an essential ingredient in the sequel. 
Let $e_{i}de_{j} \in ker \phi_{1} $.  From Eq.\ (\ref{eq:10}), we have
\[ d(e_{i}de_{j} ) = 1 \otimes e_{i} \otimes e_{j}
                   - e_{i} \otimes 1 \otimes e_{j} 
                   + e_{i} \otimes e_{j} \otimes 1. \] 
By multiplying $e_{k}$ $(k \neq i )$ to the left of $d(e_{i}de_{j}) $, 
 we obtain generators
$e_{k} \otimes e_{i} \otimes e_{j} = e_{k}de_{i} \otimes de_{j} $.
Similarly,
by multiplying $e_{k}$ $(k \neq j )$ to the right of $d(e_{i}de_{j}) $, 
 we obtain generators
 $e_{i} \otimes e_{j} \otimes e_{k} = e_{i}de_{j} \otimes de_{k} $.
Thus the remaining term $e_{i} \otimes ( 1 - e_{i} - e_{j} ) \otimes e_{j} $
 of $d(e_{i}de_{j})$, 
 which is of the form
    $e_{i}de_{1}\otimes de_{j} 
      + e_{i}de_{2}\otimes de_{j} 
      + \cdots
      + e_{i}de_{N}\otimes de_{j} $
with two terms  $e_{i}de_{i} \otimes de_{j}$ and
$e_{i}de_{j}\otimes de_{j}$
 being deleted,
is a generator of $ker \phi_{2}$.
At present, we do not have to worry about the possibility
of decomposing the generator further into simpler ones.
We note that each term of the generator begins with $e_{i}$ and 
ends with $de_{j}$. It is straightforward to proceed further in a similar 
manner.  This leads us to the following lemma.
\\

{\sc Lemma 1:} For $n \geq 2 $, $ker \phi_{n} (\neq \emptyset )$
can be generated by either single elements
$e_{i_{1}}de_{i_{2}}\otimes \cdots \otimes de_{i_{n+1}} $
or elements of the form
$
\epsilon_{I}e_{i_{1}}de_{i_{2}}\otimes \cdots \otimes de_{i_{n+1}}
+ \epsilon_{J}e_{j_{1}}de_{j_{2}}\otimes \cdots \otimes de_{j_{n+1}}
+ \cdots 
+ \epsilon_{K}e_{k_{1}}de_{k_{2}}\otimes \cdots \otimes de_{k_{n+1}} $
with $e_{i_{1}} = e_{j_{1}} = \cdots = e_{k_{1}}$ and
     $e_{i_{n+1}} = e_{j_{n+1}} = \cdots = e_{k_{n+1}}$,
where 
$\epsilon_{I}, \epsilon_{J}, \cdots ,\epsilon_{K} = \pm $.
\vspace{0.5cm}

The coefficients $\pm 1$ originate from the alternating sign 
in Eq.\ (\ref{eq:10}).
\\

{\sc Lemma 2:} Let $x \in \Omega^m $ and let $u \in ker \phi_{n} $ for
$m \geq 0, n \geq 1$. Then $x \otimes u, u \otimes x \in ker \phi_{n+m}$.
\\

{\em Proof:}  It is trivial for $ m = 0$. Let $m \neq 0$ and
let $u_{\alpha}$ be either a single element generator
 $e_{i_{1}}de_{i_{2}}\otimes \cdots \otimes de_{i_{n+1}} $
or a  generator of the form
$
\epsilon_{I}e_{i_{1}}de_{i_{2}}\otimes \cdots \otimes de_{i_{n+1}}
+ \epsilon_{J}e_{j_{1}}de_{j_{2}}\otimes \cdots \otimes de_{j_{n+1}}
+ \cdots 
+ \epsilon_{K}e_{k_{1}}de_{k_{2}}\otimes \cdots \otimes de_{k_{n+1}} $ 
 for $ker \phi_{n}$ 
($u_{\alpha} $ is a single element generator for $n = 1$).
Then for any $e_{i}de_{j} \in \Omega^1 $,
$e_{i}de_{j} \otimes u_{\alpha} =
     \delta_{ji_{1}}e_{i}du_{\alpha} \in ker \phi_{n+1}$
and
$u_{\alpha} \otimes e_{i}de_{j} =
 (-1)^{n+1}   \delta_{ii_{n+1}}(du_{\alpha})e_{j} \in ker \phi_{n+1}$.
Since $\Omega^m = (\Omega^1)^{\otimes_{A}m}$, the proof is 
completed by the induction.  QED \\

The existence of generators for $ker \phi_{n}$
with such specific forms as in Lemma 1
 gives
an $A$-bimodule which is a complement
of $ker \phi_{n} $ in $\Omega^n $ for $ n \geq 2 $.
Let $\{e_{i_{1}}de_{i_{2}}\otimes \cdots \otimes de_{i_{n+1}}
      \mid i_{k} = 1,2, \cdots , N \,\,\,(i_{k} \neq i_{k+1} ) \}$
be a basis for $\Omega^n$.
Let $S_{ij}$ be the $A$-subbimodule of $\Omega^n $ generated by
all basis elements of $\Omega^n$ with $i_{1} = i$ and $i_{n+1} = j$.
Then we have
\begin{equation}
\Omega^n = \oplus_{i,j}S_{ij}.
\label{eq:omega}
\end{equation}
By Lemma 1, we have generators of $ker \phi_{n}$, each of which belongs to
$S_{ij}$ for some $i, j$.
Let $S_{ij}^{(1)}$ be the subspace of $S_{ij}$ spanned by the 
generators of $ker \phi_{n}$ in $S_{ij}$.
Then it is clear that
\begin{equation}
ker \phi_{n} = \oplus_{ij}S_{ij}^{(1)}.
\end{equation}
Now let $S_{ij}^{(2)}$ be a complement of $S_{ij}^{(1)}$ in $S_{ij}$
and define
\begin{equation}
Q^{n} = \oplus_{ij}S_{ij}^{(2)}.
\end{equation}
Then $Q^{n}$ is not only a complement of $ker \phi_{n}$ in $\Omega^n$
but also an $A$-bimodule.
 In fact, if we let $v \in Q^n $, $v$ may be written
as $v = \sum_{ij} v_{ij}$, where $v_{ij}$ is in $S_{ij}^{(2)}$.
Let $f \in A$. Since $f(i)$
 is a complex number, both $fv = \sum_{ij}f(i)v_{ij}$ and 
$vf = \sum_{ij}f(j)v_{ij} $ belong to $Q^n$.

By construction, we have a splitting exact sequence of $A$-bimodules
\begin{equation}
0 \longrightarrow ker \phi_{1} \longrightarrow
 \Omega^1 \stackrel{\phi_{1}}{\longrightarrow} K^1 \longrightarrow 0.
\end{equation}
The splitting map $\jmath_{1} : K^1 \longrightarrow \Omega^1 $ is 
just the inclusion map.
Now we can also have a  splitting exact sequence of $A$-bimodules
for $n \geq 2$. It   is  clear that the map $\phi_{n} $ restricted to
$Q^n$
is an $A$-bimodule isomorphism.  Thus if we define a splitting map
$\jmath_{n} : K^n \longrightarrow \Omega^n $ to be the inverse
$(\phi_{n}\mid _{Q^n})^{-1}$, we have the following lemma.
\\

{\sc Lemma 3:} For $n \geq 1 $, the exact sequence of $A$-bimodules 
\[
0 \longrightarrow ker \phi_{n} \longrightarrow
 \Omega^n \stackrel{\phi_{n}}{\longrightarrow} K^n \longrightarrow 0.
\]
is split.
\\

{\em Proof: } It is enough to show that $\jmath_{n}$ 
is an $A$-bimodule homomorphism for $n \geq 2$.
Let $ \xi = e_{i_{1}}de_{i_{2}}\otimes \cdots \otimes de_{i_{n+1}} $
be a generator of $\Omega^n$.
Since $S_{i_{1}i_{n+1}} = S_{i_{1}i_{n+1}}^{(1)} \oplus
                            S_{i_{1}i_{n+1}}^{(2)}$,
 $\xi $
is expressed uniquely as $u + v $ for $u \in  S_{i_{1}i_{n+1}}^{(1)}$
and $v \in  S_{i_{1}i_{n+1}}^{(2)}$.
Thus since $\jmath_{n}( \xi + ker \phi_{n}) = v \in S^{(2)}_{i_{1}i_{n+1}} $,
  it follows  that
for any $f \in A$, 
$ \jmath_{n}(f (\xi + ker \phi_{n})) 
     = \jmath_{n}(f(i_{1})\xi + ker \phi_{n} )
      = f(i_{1})\jmath_{n}(\xi + ker \phi_{n} )
      =  f \jmath_{n}(\xi + ker \phi_{n})$. 
Similarly, $\jmath_{n}((\xi + ker \phi_{n})f) =
            \jmath_{n}(\xi + ker \phi_{n})f(i_{n+1})
           = \jmath_{n}(\xi + ker \phi_{n})f $.  
QED\\

 Let us put $ K^{0} = A $ and $\phi_{0} = \jmath_{0} =  id $.
\\

{\sc Proposition 1:}
 The following diagram

\[ \begin{array}{lllllllllllll}
A & \stackrel{d}{\rightarrow} & \Omega^{1} & \stackrel{d}{\rightarrow} &
 \cdots & \stackrel{d}{\rightarrow} & \Omega^{n-1} & \stackrel{d}{\rightarrow}
   & \Omega^{n} & \stackrel{d}{\rightarrow}
     & \Omega^{n+1} & \stackrel{d}{\rightarrow} & \cdots \\
\downarrow \phi_{0} & & \downarrow \phi_{1} & &
 & & \downarrow \phi_{n-1} 
     & & \downarrow \phi_{n}    
        & & \downarrow \phi_{n+1} & &  \\
A & \stackrel{\delta}{\rightarrow} & K^{1} & \stackrel{\delta}{\rightarrow} &
 \cdots & \stackrel{\delta}{\rightarrow}
  & K^{n-1} & \stackrel{\delta}{\rightarrow}
   & K^{n} & \stackrel{\delta}{\rightarrow}
     & K^{n+1} & \stackrel{\delta}{\rightarrow} & \cdots \\
\end{array} \]
commutes and $\delta^{2} = 0$. \\

{\em Proof:} The commutativity of the diagram is a consequence of
 the definition of $\delta $. 
For any $n = 1, 2, \cdots $,
let $u \in K^{n-1}$. Then $d\jmath_{n-1}u = \jmath_{n}v + w $
for some $v \in K^{n}$ and $w \in ker \phi_{n}$.  
Now $dw \in d(ker \phi_{n}) \subset \,\,\, ker \phi_{n+1}$, 
 and $d\jmath_{n}v = \jmath_{n+1}v_{1} + v_{2} $ for some $v_{1} \in K^{n+1} $ 
   and $v_{2} \in  \,  ker \phi_{n+1}$.
Thus
  from $d^{2}\jmath_{n-1}u = 0 $, it follows that 
$ \jmath_{n+1}v_{1} + v_{2} + dw = 0$.
But $\jmath_{n+1}v_{1} = -v_{2} - dw = 0$ since 
$\jmath_{n+1}K^{n+1} \cap ker \phi_{n+1} = \{ 0 \}$. 
  Hence $\delta^{2}u = v_{1} = 0$. QED \\

Let us define $P_{n,m} : K^n \otimes _{A} K^{m} \longrightarrow K^{n+m} $,
$\omega \otimes \omega' \mapsto \omega \cdot \omega' $, to be
$P_{n,m} \equiv \phi_{n+m} \circ (\jmath_{n} \otimes \jmath_{m} ) $
for $n, m \geq 0 $.
Note that the multiplication in $\Omega(A)$ is implicitly incorporated
from Eq.\ (\ref{eq:15}).
The map $P_{n,m}$ is well-defined.  In fact, 
it is trivial for the cases where $n = 0$ and  $ m = 0$. Otherwise,
let us take any splitting maps
$\tilde{\jmath_{n}}$ and $\tilde{\jmath_{m}}$. Then 
$\jmath_{n} \omega \otimes \jmath_{m} \omega'
 - \tilde{\jmath_{n}} \omega \otimes \tilde{\jmath_{m}}\omega'
=  \jmath_{n}\omega \otimes ( \jmath_{m} \omega' - \tilde{\jmath_{m}}\omega')
 +  ( \jmath_{n} \omega - \tilde{\jmath_{n}}\omega ) \otimes 
          \tilde{\jmath_{m}}\omega' $.
Since $ \jmath_{m} \omega' - \tilde{\jmath_{m}}\omega' \in ker \phi_{m}$
and $\jmath_{n}\omega \in \Omega^n $ etc., 
$\jmath_{n} \omega \otimes \jmath_{m} \omega'
 - \tilde{\jmath_{n}} \omega \otimes \tilde{\jmath_{m}}\omega'$ 
belongs to $ker \phi_{n+m} $ by  Lemma 2.
Thus $\phi_{n+m} (\jmath_{n} \omega \otimes \jmath_{m} \omega')
= \phi_{n+m}( \tilde{\jmath_{n}} \omega \otimes \tilde{\jmath_{m}}\omega')$.
Also, by the $A$-bilinearity of the splitting maps, 
$P_{n,m}(\omega f \otimes \omega' )
 = P_{n,m}(\omega \otimes f \omega' ) $ 
for any $f \in A$.
\\

Moreover, $P_{n,m} $ is associative.
\\

{\sc Lemma 4:} Let $\omega \in K^n, \omega' \in K^m $ and $\omega'' \in K^l $
for $n, m, l \geq 0$.
Then
\[ (\omega \cdot \omega')\cdot \omega''
     = \omega \cdot (\omega'\cdot \omega'').
\]

{\em Proof:} For $n = m = l = 0$, it is just the associativity of $A$. 
 Otherwise, we have
\begin{eqnarray*}
(\omega \cdot \omega')\cdot \omega''
  & = & P_{n+m,l}(( \omega \cdot \omega')\otimes \omega'') \\
  & = & \phi_{n+m+l}[\jmath_{n+m}(\phi_{n+m}(\jmath_{n}\omega
         \otimes \jmath_{m}\omega'))\otimes \jmath_{_{l}}\omega'' ] \\
  & = & \phi_{n+m+l}[\jmath_{n}\omega \otimes \jmath_{m+l}(\phi_{m+l}
        (\jmath_{m}\omega'
         \otimes \jmath_{_{l}}\omega'')) ] \\
 & = & P_{n,m+l}(( \omega \otimes (\omega' \cdot \omega'')) =
          \omega \cdot (\omega' \cdot \omega''). 
\end{eqnarray*}
The third equality comes from the fact that
\[ \jmath_{n+m}(\phi_{n+m}(\jmath_{n}\omega
         \otimes \jmath_{m}\omega'))\otimes \jmath_{_{l}}\omega''
- \jmath_{n}\omega \otimes \jmath_{m+l}(\phi_{m+l}
        (\jmath_{m}\omega'
         \otimes \jmath_{_{l}}\omega''))\]
\[= [\jmath_{n+m}(\phi_{n+m}(\jmath_{n}\omega
         \otimes \jmath_{m}\omega'))-\jmath_{n}\omega
             \otimes \jmath_{m}\omega']\otimes \jmath_{_{l}}\omega''
- \jmath_{n}\omega \otimes [\jmath_{m+l}(\phi_{m+l}
        (\jmath_{m}\omega'
         \otimes \jmath_{_{l}}\omega''))
 - \jmath_{m}\omega' \otimes \jmath_{_{l}}\omega''], 
\]
which is in $ker \phi_{n+m+l} $ by  Lemma 2 since
$\jmath_{n+m}(\phi_{n+m}(\jmath_{n}\omega
         \otimes \jmath_{m}\omega'))-\jmath_{n}\omega
             \otimes \jmath_{m}\omega' \in  ker \phi_{n+m}$
and $ \jmath_{_{l}}\omega'' \in \Omega^l $ etc..
 QED \\

{\sc Lemma 5:} $\phi_{n+m} = P_{n,m}(\phi_{n} \otimes \phi_{m})$
for $n, m \geq 0 $.
\\

{\em Proof:} It is trivial for $n = m = 0$. Otherwise,
let $\alpha \otimes \beta \in \Omega^n \otimes_{A}
\Omega^m = \Omega^{n+m} $.
\begin{eqnarray*}
P_{n,m} (\phi_{n} \otimes \phi_{m} )(\alpha \otimes \beta )
 & = & (\phi_{n+m} \circ (\jmath_{n} \otimes \jmath_{m}))(\phi_{n}\alpha
      \otimes \phi_{m}\beta ) \\
 & = & \phi_{n+m}(\jmath_{n}\phi_{n}\alpha \otimes \jmath_{m}\phi_{m}\beta ).
\end{eqnarray*}
Now $ \jmath_{n}\phi_{n}\alpha \otimes \jmath_{m}\phi_{m}\beta 
 - \alpha \otimes \beta 
= \jmath_{n}\phi_{n}\alpha \otimes (\jmath_{m}\phi_{m}\beta 
  - \beta ) + 
 (\jmath_{n}\phi_{n}\alpha - \alpha ) \otimes \beta  $,
which is in $ker \phi_{n+m} $ by  Lemma 2. 
Hence
$\phi_{n+m}(\jmath_{n}\phi_{n}\alpha \otimes \jmath_{m}\phi_{m}\beta )
 = \phi_{n+m}(\alpha \otimes \beta )$.
QED\\

Now we are ready to show the fact that
the operator $\delta $ satisfies the Leibniz rule as $d$ in 
Eq.(\ref{eq:2}).\\

{\sc Proposition 2:} Let $\omega \in K^n $ and $\omega' \in K^m$
for $n, m \geq  0 $.  Then
\[
 \delta(\omega \cdot \omega ') = (\delta\omega)\cdot \omega ' +
 (-1)^{n}\omega \cdot\delta\omega '.  
\] 

 {\em Proof:} For  $n, m \geq 0 $, it follows that
\begin{eqnarray*}
\delta (\omega \cdot \omega' ) 
   & = & \delta P_{n,m}(\omega \otimes \omega')
         = \delta \phi_{n+m} (\jmath_{n} \otimes \jmath_{m})
           (\omega \otimes \omega')\\
  &  = & \phi_{n+m+1}d(\jmath_{n}\omega \otimes \jmath_{m}\omega' )
    = \phi_{n+m+1}(d\jmath_{n}\omega \otimes \jmath_{m}\omega' + (-1)^{n}
            \jmath_{n}\omega \otimes d\jmath_{m} \omega'  )\\
   & = & P_{n+1,m}(\phi_{n+1}d\jmath_{n}\omega \otimes \phi_{m}\jmath_{m}\omega' )
        + (-1)^{n} P_{n,m+1}
            (\phi_{n}\jmath_{n}\omega \otimes \phi_{m+1}d\jmath_{m}\omega' )\\
   & =  &  (\delta\omega)\cdot \omega ' +
              (-1)^{n}\omega \cdot\delta \omega ',
\end{eqnarray*}
where the fifth equality comes from Lemma 5. QED
\vspace{0.5cm}

Thus we have proved the following.
\vspace{0.5cm}

{\sc Proposition 3:}
$ K(A) = A \oplus K^{1} \oplus K^{2} \oplus \cdots $
 is  a differential algebra.
\vspace{0.5cm}

It is well known$^{6}$  that $ <d\Omega^{n}> = \Omega^{n+1}$.
  The same relation
 holds for $ \delta $.  In fact,
 $K^{n+1} \supset \phi_{n+1} (<d\jmath_{n}K^{n}>) $.
To show  
$ K^{n+1} \subset \,\,\, \phi_{n+1}(<d\jmath_{n}K^{n}>)$,
we first observe that any element of $\Omega^{n+1}$ is 
expressed ( not necessarily uniquely ) as a 
sum of an element in $<d\,ker \phi_{n}>$
and an element in $< d\,\jmath_{n}K^{n}>$ since 
$\Omega^{n+1} = < d\Omega^{n} > $.
Thus if we let
 $v \in K^{n+1} $,  we have $\jmath_{n+1}(v) = v_{1} + v_{2}$
for some $v_{1} \in \, <d\,ker \phi_{n}>$ 
  and $v_{2} \in \, <d\,\jmath_{n}K^{n}>$.
Then $v = \phi_{n+1}\jmath_{n+1}(v) = \phi_{n+1}(v_{1} + v_{2})
        =\phi_{n+1}(v_{2}) \in \phi_{n+1}(<d\,\jmath_{n}K^{n}>)$.
Thus
 $K^{n+1} = \phi_{n+1} (<d\jmath_{n}K^{n}>) $.
On the other hand, $\phi_{n+1}(<d\jmath_{n}K^n >)
 = <\phi_{n+1}d\jmath_{n}K^n > = < \delta K^n > $ since
$\phi_{n+1} $ is an $A$-bimodule homomorphism.  Thus we have   
 $ <\delta K^{n}> = K^{n+1} $. \\

From now on, we shall write $e_{i_{1}i_{2}\cdots i_{n}} $ for
$e_{i_{1}}de_{i_{2}} \otimes  \cdots \otimes de_{i_{n}} $
for simplicity, according to Ref.\ 6, with the convention that
$e_{i_{1}i_{2}\cdots i_{n}} $ = $0$ if $i_{k} = i_{k+1}$ for some $k$.
Also we   shall  often write $v $
 for $ x \,+ ker \phi_{n}  $ of $K^n$ if 
 $\jmath_{n}(x \,+ ker \phi_{n}) = v $
for a given splitting map $\jmath_{n}$. 
\\

{\sc Example 1:} 
 Let $\cal V $ be a set $\{x_{1}, x_{2}, x_{3} \}$ of three points.
          If we let 
$K^{1} $ be the space generated by $\{ e_{12}, e_{23}, e_{13} \}$,
then  $ker \phi _{1}\,=\,\, < \{ e_{21}, e_{32}, e_{31} \}> $. From this, 
it follows that 
 $ker \phi_{2}  $ is generated by $  \Omega ^{2} - <\{e_{123} \}> $,
and $K^{2} =\,\,<e_{123}> $ with the obvious splitting map $\jmath_{2}$.
 Moreover, $K^{n} = 0 $ for $n \geq  3 $. 

         Similarly, if we  let 
$K^{1} $ be the space generated by
 $\{ e_{12}, e_{21}, e_{23}, e_{32} \} $,
then  $ker \phi_{1} = \,\,< \{ e_{13}, e_{31} \}> $ and
 $ker \phi_{2}  $ is generated by $ \Omega ^{2} - 
<\{e_{121}, e_{212}, e_{232}, e_{323} \}> $.   
Also $K^{2} =\,\, <\{e_{121}, e_{212}, e_{232}, e_{323} \}> $, etc..    
 \\
  
{\sc Example 2:} 
Let  $\cal V $ be a set  $\{x_{1}, x_{2}, x_{3}, x_{4} \}$ of
four points
and let 
$K^{1} $ be the space generated by $\Omega^{1} -
 <\{ e_{14}, e_{41} \}> $.
Then  $ker \phi_{1} \,\,= \,\,< \{ e_{14}, e_{41} \}> $ and
 $ker \phi_{2}  $ is generated by the following 12 elements 
\[
e_{141}, e_{142}, e_{143}, e_{214}, e_{241}, \\    
e_{314}, e_{341}, e_{412}, e_{413}, e_{414}, \\
e_{124} + e_{134}, e_{421} + e_{431}.           
\]
Note that the two elements
$ e_{124} + ker \phi_{2}  $ and $e_{134} + ker \phi_{2} $
 are not linearly independent in $K^2$
since $e_{124} + e_{134}$ belongs to $ker \phi_{2} $. 
Rather, $e_{124} = -e_{134} $ modulo $ ker \phi_{2} $.\\

In the case of the universal differential algebra $\Omega (A)$, 
it is known that the sequence
\begin{equation}
A \stackrel{d}{\rightarrow}
  \Omega^{1} \stackrel{d}{\rightarrow}
    \Omega^{2} \stackrel{d}{\rightarrow}
      \cdots \stackrel{d}{\rightarrow}
          \Omega^{n} \stackrel{d}{\rightarrow}
           \cdots.
\label{eq:17}       
\end{equation}
is exact$^{6}$. However, this is not true for $\delta $.
The counter-example can be seen in the following.
\vspace{0.5cm}

{\sc Example 3:}
 Let $\cal V $ be a set $\{x_{1}, x_{2}, x_{3} \}$ of three points
and let $K^{1} $ be the space generated by $\{ e_{12}, e_{23}, e_{31}\} $.
Then $ker \phi_{1} = < \{ e_{21}, e_{32}, e_{13} \}> $.
It is straightforward to see that 
$ker \phi_{2}  = \Omega ^{2} $.   
Hence $ K^{n} = 0 $ 
for $ n \geq 2 $.
  We note that there is an element, say 
$e_{12}\otimes e_{23} = e_{123}$, which  is in $K^1 \otimes_{A}  K^1$, but
 not in $K^{2}$.

  From the fact that the dimensions  of both $A$  and $K^1$ 
are the same in
the  sequence  
$ A \stackrel{\delta}{\rightarrow}
   K^{1} \stackrel{\delta}{\rightarrow}
    K^{2} = 0 
$, the map $\delta : A \longrightarrow K^1 $ is   
 an isomorphism if the sequence is exact.  But this contradicts the fact that $\delta 1 = 0$.
Thus the sequence is not exact.
\vspace{1cm}

{\bf \Large III. LINEAR CONNECTIONS ON GRAPHS}
\vspace{0.5cm}

{\bf \large A. Linear connections}
\vspace{0.5cm}

Let $A$ be an associative algebra and  $E$ be
an $A$-bimodule. One may impose a reality condition once  
 $A$ is given as a $*$-algebra$^{12}$. 
 Let $(\Omega^{*}, d ) $ 
 be a
 differential calculus over $A$. If  $\Omega^1$ 
is an $A$-bimodule,  we can define a left and a right connection 
on $E$.
A left connection on $E$ is defined to be
 a linear map
 \begin{equation}
D : E \longrightarrow  \Omega^1 \otimes_{A} E  
\label{eq:311}
\end{equation}
satisfying
\begin{equation}
D ( f\omega ) = df \otimes \omega + f D \omega 
\label{eq:312}
\end{equation}
 for any $f \in A $  and $ \omega \in E $.
One can also define a right connection on $E$ to be a linear map
 \begin{equation}
D : E \longrightarrow  E \otimes_{A} \Omega^1  
\label{eq:313}
\end{equation}
satisfying 
\begin{equation}
D ( \omega f ) = ( D \omega)f +  \omega \otimes df . 
\label{eq:314}
\end{equation}

In Ref.\ 15, a definition of a bimodule connection is proposed.
A bimodule connection  on $E$ is a left connection $D $   
such that for any $f \in A $  and $\omega \in E $, $D (\omega f)$ 
is of the form
\begin{equation}
D ( \omega f ) = ( D \omega)f + \sigma( \omega \otimes df), 
\label{eq:315}
\end{equation}
where $\sigma $ is a map from $ E \otimes_{A} \Omega^{1} $
to $\Omega^1 \otimes_{A} E $ generalizing the permutation.
In particular, if we take $E = \Omega^1 $, the bimodule connection
is called a linear connection,  which we are mainly 
concerned with in this work.  

Let $\Omega^2 $ be the $A$-bimodule of two-forms and 
let $\pi : \Omega^1 \otimes_{A} \Omega^1 \longrightarrow \Omega^2 $ be 
a linear map satisfying $\pi (\omega \otimes \omega' ) 
= \omega \cdot \omega' $ where $\cdot $ is the multiplication map 
between forms. 
For the consistencey of the definition of a linear connection, the map
$\sigma : \Omega^1 \otimes _{A}  \Omega^1
 \longrightarrow  \Omega^1 \otimes _{A} \Omega^1 $
is assumed to be $A$-bilinear, i.e. 
for $f \in A$ and $\omega, \omega' \in \Omega^1$,
\begin{equation}
\sigma (f \omega \otimes \omega' ) = 
f \sigma ( \omega \otimes \omega' ),
\hspace{1cm}
\sigma ( \omega \otimes \omega'f ) = 
 \sigma ( \omega \otimes \omega' ) f.  
      \label{eq:415}
\end{equation}
Moreover, $\sigma $ is assumed to satisfy the following
\begin{equation}
\pi \circ (\sigma + 1 ) = 0.
      \label{eq:416}
\end{equation}
The relation in Eq.\ (\ref{eq:416}) is a necessary and sufficient condition for
the torsion $T$ of the connection $D$,
 defined by $ T = d - \pi \circ D $, to be 
  $A$-bilinear ( more precisely, right $A$-linear).

A linear connection $D $ can be extended to 
two linear maps $D_{1} : \Omega^1 \otimes_{A} \Omega^1 \longrightarrow 
\Omega^2 \otimes_{A} \Omega^1,$ and $
 D_{2} :\Omega^1 \otimes_{A} \Omega^1 \longrightarrow 
\Omega^1 \otimes_{A} \Omega^1 \otimes_{A} \Omega^1 $, respectively, satisfying
\begin{equation} 
D_{1}(\omega \otimes \omega' )  =  d \omega \otimes \omega'
   - \pi_{12}( \omega \otimes D \omega' )  
\label{eq:417}
\end{equation} 
and
\begin{equation}  
D_{2}(\omega \otimes \omega' )  =  D \omega \otimes \omega'
   + \sigma_{12}( \omega \otimes D \omega')  
\label{eq:418}
\end{equation} 
for $\omega, \omega' \in \Omega^1 $,
where $\pi_{12} = \pi \otimes 1 $ and 
$\sigma_{12} = \sigma \otimes 1 $.

It is easy to see that the linear map $D_{1} \circ D $ 
 is a left $A$-linear. Also 
$\pi_{12} D_{2} \circ D $ is a left $A$-linear  if
the torsion $T = 0$.
For the right $A$-linearity
of  $D_{1} \circ D $ 
 and
$\pi_{12} D_{2} \circ D $,
 there is not yet  a widely accepted prescription 
 even though it seems to be an
essential property for the concept of a curvature.
However, one prescription to obtain an $A$-bilinear curvature
has been proposed recently$^{19}$, which we shall use in this work.
Especially, since $\Omega^1 $ is free for graphs, one can construct
the curvature invariants from the linear connection$^{19}$. 
In the next subsections, we shall 
calculate  linear connections 
and  curvatures explicitly on graphs 
with respect to the natural basis
for the general nonzero
torsion case. 
There are some other models$^{20-23}$
 for which linear connections and curvatures 
are calculated mostly without torsion. 
\vspace{0.5cm}

{\bf \large B. Linear connections on complete symmetric graphs }
\vspace{0.5cm}

Let $A$ be the associative algebra of complex functions
 on a set ${\cal V}$ of $N$ points and let $\Omega (A)$ be the 
universal differential algebra,
 which corresponds to a complete symmetric graph,
 introduced in Sec.\ II.\ B.
Since $\pi = 1 $ for the universal differential algebra $\Omega(A)$,
we take $\sigma = -1 $ from Eq.\ (\ref{eq:416})
 for the $A$-bilinearity of the torsion $T$.
Thus 
 a linear connection is given by 
 a linear map $D : \Omega^1
 \longrightarrow \Omega^1 \otimes_{A} \Omega^1 $ 
 satisfying
\begin{equation}
 D ( f\omega ) = df \otimes \omega + f D \omega,  
\hspace{1cm}
D ( \omega f ) = (D \omega )  f - \omega \otimes df 
\end{equation} 
for any $f \in A $ and $\omega \in \Omega^1 $.
Moreover, $D_{1}$ and $D_{2}$ satisfy
\begin{equation}
D_{1}(\omega \otimes \omega') = d\omega \otimes \omega'
- \omega \otimes D\omega',
\hspace{0.5cm}
D_{2}(\omega \otimes \omega') = D\omega \otimes \omega'
- \omega \otimes D\omega'.
\end{equation}
In this case, we have $D_{1} = D_{2} + T \otimes 1$.

Since $ \Omega^1 \otimes _{A} \Omega^1 = \Omega^2 $ is
 a vector space with a  basis 
$\{ e_{ijk}  \mid i, j, k = 1, 2, \cdots , N \,\,\,
 (i \ne j , j \neq k ) \} $,
we may put 
\begin{equation}
 D(e_{ij}) = \sum_{k,l,m}\Gamma_{ij}^{klm}e_{klm} 
\label{eq:321}
\end{equation}
for some numbers $\Gamma_{ij}^{klm} $'s. 
Here let us put $\Gamma_{ij}^{klm} = 0 $
if $ i = j $ or $ k = l$ or $ l = m $. 

Now from the following two expressions of $D(e_{ij})$
\begin{eqnarray}
D(e_{ij}) & = & D(e_{i}de_{j}) 
                =  de_{i} \otimes de_{j} + e_{i}D(de_{j})  \nonumber \\
        &=& \sum_{a}(e_{aij} - e_{iaj} + e_{ija}) 
               + \sum_{m,b,c}(\Gamma_{mj}^{ibc} - \Gamma_{jm}^{ibc})e_{ibc} 
\end{eqnarray}
and
\begin{eqnarray}
D(e_{ij}) & = & D(e_{i}de_{j}) = -D(de_{i}e_{j}) 
                = - D(de_{i})e_{j} + de_{i} \otimes de_{j}   \nonumber \\
        &=& \sum_{a}(e_{aij} - e_{iaj} + e_{ija}) 
               - \sum_{m,a,b}(\Gamma_{mi}^{abj} - \Gamma_{im}^{abj})e_{abj}, 
   \label{eq:324}
\end{eqnarray}
we have $\Gamma_{ij}^{abc} = 0 $ except for
\begin{eqnarray}
\Gamma_{ij}^{ija} &=& 1 \,\,\,\,\, (a \ne j ), \nonumber \\
\Gamma_{ij}^{aij} &=& 1 \,\,\,\,\, (a \ne i ) 
\label{eq:325}
\end{eqnarray}
with $\Gamma_{ij}^{iaj} ( a \ne i,j ) $ undetermined.
From these values of $\Gamma_{ij}^{klm}$'s, we obtain the following lemma. 
\vspace{0.5cm}

{\sc Lemma 6:} For any $i, j \,\, (i \ne j )$, $D(e_{ij})$ and $T(e_{ij})$ are 
of the following forms
\begin{eqnarray*}
  D(e_{ij}) &=& de_{i} \otimes de_{j} + 
               \sum_{a}( 1 + \Gamma_{ij}^{iaj} )e_{iaj},\\
  T(e_{ij}) & = & - \sum_{a}(1 + \Gamma_{ij}^{iaj})e_{iaj}.
\end{eqnarray*}
\vspace{0.5cm}

Thus it is obvious that the torsion $T = d - D $ is 0 if and only if 
$\Gamma_{ij}^{iaj} = -1 $ for all $i,j ( i \neq j )$ and $ a (\neq i, j)$.
Moreover, 
if the torsion $T$ is 0, the curvatures
 $D_{1} \circ D = D_{2} \circ D = 0$ in
a complete symmetric graph, which is already known (see, e.g., Ref.\ 23).
However, we shall calculate  curvatures 
for the nonvanishing torsion $T$. 

We can have a  general form of the curvature $D_{1}\circ D $
for any linear connection $D$ obeying Eq.(\ref{eq:321}):
\begin{eqnarray}
D_{1} D e_{ij} &=& - \sum_{k,l,m}D_{1}
                 (\Gamma_{ij}^{klm}de_{k}\otimes e_{l}de_{m})
                             \nonumber \\
        &=&\sum_{k,l,m}  
              \Gamma_{ij}^{klm}de_{k}\otimes D(e_{l}de_{m}) \nonumber \\
        &=& \sum_{a,b,c,d}
           \sum_{l,m}(\Gamma_{ij}^{blm}\Gamma_{lm}^{bcd} 
           - \Gamma_{ij}^{alm}\Gamma_{lm}^{bcd})e_{abcd} \\
        &=& \sum_{a,b,c,d} 
           \Omega_{ij}^{abcd}e_{abcd}, \nonumber
\label{eq:328}
\end{eqnarray}
where 
\begin{equation}
\Omega_{ij}^{abcd} = \sum_{l,m}
  ( \Gamma_{ij}^{blm}\Gamma_{lm}^{bcd} - \Gamma_{ij}^{alm}\Gamma_{lm}^{bcd}).
   \label{eq:329}
\end{equation}

Using $\Gamma_{ij}^{klm}  = \delta_{i}^{k}\delta_{j}^{l} + 
                            \delta_{i}^{l}\delta_{j}^{m} +
                            \Gamma_{ij}^{klm}\delta_{i}^{k}\delta_{j}^{m} $
in Eq.\ (\ref{eq:325}),  we obtain $D_{1}De_{ij}$.
Also, $D_{2}De_{ij}$ can be obtained from Lemma 6 and the
following observations 
\begin{equation}
D(de_{j}) = \sum_{a} \sum_{m \ne j}
       [(1+ \Gamma_{mj}^{maj})e_{maj} - (1+ \Gamma_{jm}^{jam})e_{jam}],
\label{eq:330}
\end{equation}
and for $ a \ne i, j $,
\begin{equation}
de_{i} \otimes de_{a} \otimes de_{j} = 
 \sum_{l}(e_{liaj} - e_{ilaj} + e_{ialj} - e_{iajl}).
\label{eq:331}
\end{equation}
The results are as follows.
\vspace{0.5cm}

{\sc Proposition 4:}  
\begin{eqnarray*}
D_{1}De_{ij} & = &  - \sum_{l \neq i}( 1 + \Gamma_{ij}^{ilj} )e_{ijlj}
                 +  \sum_{l\neq j} (1 - \Gamma_{ij}^{ilj}\Gamma_{lj}^{lij})
                        e_{ilij} \\
             & & - \sum_{l \neq j}\sum_{m \neq i} 
                   (\Gamma_{ij}^{imj} + \Gamma_{ij}^{ilj}\Gamma_{lj}^{lmj})     
                  e_{ilmj} \\
             & & - \sum_{l}\sum_{m \neq j}
                   ( 1 + \Gamma_{jm}^{jlm})e_{ijlm}
                 -  \sum_{l}\sum_{m} 
                   ( 1 + \Gamma_{ij}^{ilj})e_{iljm},
\end{eqnarray*}
and 
\begin{eqnarray*}
D_{2}De_{ij} & = &  - \sum_{l \neq i}( 1 - \Gamma_{ij}^{ilj}
                         \Gamma_{il}^{ijl} )e_{ijlj}
                 +  \sum_{l\neq j} (1 - \Gamma_{ij}^{ilj}\Gamma_{lj}^{lij})
                        e_{ilij} \\
             & & + \sum_{l \neq j}\sum_{m \neq i} 
                   (\Gamma_{ij}^{imj}\Gamma_{im}^{ilm}
                    - \Gamma_{ij}^{ilj}\Gamma_{lj}^{lmj})     
                  e_{ilmj} \\
             & & - \sum_{l}\sum_{m \neq j}
                   ( 1 + \Gamma_{jm}^{jlm})e_{ijlm}
                 +  \sum_{l\neq i}\sum_{m} 
                   ( 1 + \Gamma_{li}^{lmi})e_{lmij}.
\end{eqnarray*}

If $D_{1} \circ D $ is $A$-bilinear, the torsion $T$ should vanish, i.e.
all $\Gamma_{jm}^{jlm} = -1 $
from the vanishment of the $e_{ijlm}$
or $e_{iljm} $ term in $D_{1}De_{ij}$. Similarly for $D_{2}\circ D$.
 We thus have
\vspace{0.5cm}

{\sc Corollary 1:} Let the torsion $T$ be $A$-bilinear 
                 on a complete symmetric graph.  Then  
the necessary and sufficient condition for
 $D_{1}\circ D $ to be $A$-bilinear is $T$ = 0.
Similarly for $D_{2}\circ D$.
\vspace{0.5cm}

One prescription to get an $A$-bilinear curvature for the nonzero torsion $T$
is to factor out all those elements that do not satisfy the desired
condition.  We refer to Ref.\ 19 for a recent discussion about
 this prescription.
Thus if we factor out the terms $e_{ijlm}$ and $e_{iljm}$
 of $D_{1}De_{ij}$, 
we have an $A$-bilinear curvature, denoted by $Curv_{1} $, since
the remaining terms belong to the $A$-bimodule $S_{ij}$.
Similarly, we have an $A$-bilinear curvature denoted by $Curv_{2} $
from $\pi_{12} D_{2}De_{ij} = D_{2}De_{ij}$.
\vspace{0.5cm}

{\sc Corollary 2:} For all $i, j \,\,\, (i \neq j),$
\begin{eqnarray*}
 Curv_{1}(e_{ij}) & = &  - \sum_{l \neq i}( 1 + \Gamma_{ij}^{ilj} )e_{ijlj} 
               +  \sum_{l\neq j} (1 - \Gamma_{ij}^{ilj}\Gamma_{lj}^{lij})
                        e_{ilij} \\
             & & - \sum_{l \neq j}\sum_{m \neq i} 
                   (\Gamma_{ij}^{imj} + \Gamma_{ij}^{ilj}\Gamma_{lj}^{lmj})     
                  e_{ilmj},  
\end{eqnarray*}
and
\begin{eqnarray*}
 Curv_{2}(e_{ij})  & = &  - \sum_{l \neq i}( 1 - \Gamma_{ij}^{ilj}
                         \Gamma_{il}^{ijl} )e_{ijlj}
                 +  \sum_{l\neq j} (1 - \Gamma_{ij}^{ilj}\Gamma_{lj}^{lij})
                        e_{ilij} \\
             & & + \sum_{l \neq j}\sum_{m \neq i} 
                   (\Gamma_{ij}^{imj}\Gamma_{im}^{ilm}
                    - \Gamma_{ij}^{ilj}\Gamma_{lj}^{lmj})     
                  e_{ilmj}.
\end{eqnarray*}

Now let us consider  an interesting special case where
\begin{equation}
\Gamma_{ij}^{klm} =
 \Gamma_{\varrho(i)\varrho(j)}^{\varrho(k)\varrho(l)\varrho(m)}
\label{physical}
\end{equation}    
for any permutation  $\varrho $  on
the set $\cal V $ of $N$ points.
This is motivated by the fact that a complete symmetric graph 
$(A, \Omega^1 )$
is invariant under the permutations $\varrho $ on $\cal V $
in the sense that
 the  transformed bases $\{e_{\varrho(i)}\}$    for $A$ and 
$\{ e_{\varrho(i)\varrho(j)} \}$
for $\Omega^1 $ are  equivalent respectively to the original ones since 
the graph $(A, \Omega^1 )$ is complete and symmetric.
In this special case,
 let us put $1 + \Gamma_{ij}^{iaj} = \gamma $   
for all $i, j \,\,\, (i \neq j)$ and $a (\neq i, j) $.
Then the curvatures in the above can be written immediately  as follows.
\vspace{0.5cm}

{\sc Corollary 3:} For all $i, j \,\,\,(i \neq j)$,
\begin{eqnarray*}
D_{1} D e_{ij}
        &=& \gamma \{ -\sum_{l\ne i}e_{ijlj}
                       + (2-\gamma)\sum_{l\ne j}e_{ilij}
                       + (1-\gamma)\sum_{l\ne j}\sum_{m \ne i}e_{ilmj}
                              \\
       & & - \sum_{l}\sum_{m \ne j}e_{ijlm}
                     - \sum_{l}\sum_{m }e_{iljm} \}, \\
Curv_{1}(e_{ij}) & = & \gamma \{  -\sum_{l\ne i}e_{ijlj}
                       + (2-\gamma )\sum_{l\ne j}e_{ilij} 
                     + (1-\gamma)\sum_{l\ne j}\sum_{m \ne i}e_{ilmj}
                        \}, 
\end{eqnarray*}
and
\begin{eqnarray*}
D_{2} D e_{ij}
       & = & \gamma \{ (\gamma - 2)\sum_{l\ne i}e_{ijlj}
                       + (2-\gamma)\sum_{l\ne j}e_{ilij} 
                   -  \sum_{l}\sum_{m \ne j}e_{ijlm}  \\
                &&     +  \sum_{l \neq i}\sum_{m }e_{lmij} \}, \\
 Curv_{2}(e_{ij}) &= & \gamma (\gamma - 2) \{\sum_{l\ne i}e_{ijlj}
                       - \sum_{l\ne j}e_{ilij} \}. 
\end{eqnarray*}
\vspace{0.5cm}

 The parameter $\gamma$ determines a connection and hence 
 curvatures. These one-parameter families of connections and curvatures 
on a graph are closely analogous to 
those on matrix geometries$^{20}$ and
those on the ordinary quantum plane$^{23}$.
However, the torsion is not zero in general in this work. In fact,
the torsion also depends on the papameter
since $Te_{ij} = - \gamma \sum_{a}e_{iaj} $ from Lemma 6. 
A surprising result that $Curv_{1} = Curv_{2} $ arises
when $\gamma = 1 $ for  complete symmetric graphs.   

We define a metric $g : \Omega^1 \otimes _{A} \Omega^1 \longrightarrow A$
on $(A, \Omega^1 )$ to be an $A$-bilinear nondegenerate map.
By nondegenerate, we mean $ g(\omega \otimes \omega' ) = 0 $ 
for all $\omega \in \Omega^1 $ implies $\omega' = 0$ and
 $ g(\omega \otimes \omega' ) = 0 $ 
for all $\omega' \in \Omega^1 $ implies $\omega = 0$.
Then we have $g(e_{ij} \otimes e_{jk} ) = \mu_{i} e_{i} \delta_{ik} $
for some constant $\mu_{i} $.
From the same motivation as that for Eq.\ (\ref{physical}), we assume that
$\mu_{i} $'s are the same, say $\mu $.
Now 
 if we define a metric-compatible connection
to be a linear connection satisfying 
$d \circ g = (1 \otimes g ) \circ D_{2}$ as usual$^{20-23}$, 
then it is straightforward
to see that $dg(e_{iji}) = \mu de_{i} $ and
$( 1 \otimes g )D_{2}(e_{iji}) = \mu (de_{i}e_{i} - e_{i}de_{j})$.
Hence
 there is no metric-compatible connection in general except
the $N = 2 $ case.
This fact tells us that a metric in this sense is not so useful
for graphs.

From the connections and curvatures of a complete symmetric graph,
 we can induce
those of subgraphs. We shall do this in the next subsection.
\vspace{0.5cm}

{\bf \large C. Linear connections on graphs}
\vspace{0.5cm}

Let $(A, K^1 ) $ be a subgraph of $(A, \Omega^1 )$.  
Now we define a linear connection $\nabla : K^1 \longrightarrow
K^1 \otimes_{A} K^1 $ on a graph $(A, K^1 ) $ by the composite map
$\nabla \equiv (\phi_{1} \otimes \phi_{1} ) \circ D \circ \jmath_{1} $
where $\jmath_{1}  : K^1 \longrightarrow \Omega^1 $ is the splitting map 
 and $\phi_{1} \otimes \phi_{1} : \Omega^1 \otimes_{A}  
\Omega^1 \longrightarrow K^1 \otimes_{A} K^1 $
is the projection map.  
We also define a generalized permutation 
$\tau : K^1 \otimes_{A} K^1 \longrightarrow K^1 \otimes_{A} K^1 $
by 
$ \tau \equiv (\phi_{1}\otimes \phi_{1} )\circ \sigma \circ
(\jmath_{1} \otimes \jmath_{1})$.
Let $\triangle \equiv \Omega^1 \otimes_{A} ker \phi_{1} + ker \phi_{1}
\otimes_{A} \Omega^1 $.
Then we have
\vspace{0.5cm}

{\sc Proposition 5: } For $f \in A $ and $\omega \in K^1 $, 

(1) $ \nabla (f\omega ) = \delta f \otimes \omega + f\nabla \omega $ and

(2) $\nabla (\omega f )    = (\nabla \omega )f + \tau 
                              (\omega \otimes \delta f) $
        if $\sigma $ preserves $\triangle$, i.e.
 $\sigma(\triangle ) \subset \triangle $. 
\vspace{0.5cm}

{\em Proof:} (1) $\nabla (f \omega ) =
 (\phi_{1} \otimes  \phi_{1} )
      (df \otimes \jmath_{1} \omega + f D\jmath_{1}\omega )
= \phi_{1}df \otimes \phi_{1}(\jmath_{1}\omega) + f(\phi_{1}\otimes \phi_{1})
(D\jmath_{1}\omega )
= \delta f \otimes \omega + f\nabla \omega $.
(2) $\nabla(\omega f ) = (\phi_{1} \otimes \phi_{1} )
((D\jmath_{1}\omega )f + \sigma( \jmath_{1}\omega \otimes df )) 
       = (\nabla \omega )f 
+ (\phi_{1} \otimes \phi_{1} )\sigma (\jmath_{1}\omega \otimes
\jmath_{1}\delta f + \jmath_{1}\omega \otimes u )
 $ for some $u \in ker \phi_{1} $
since $\phi_{1}( df - \jmath_{1}\delta f ) = 0$. 
By assumption, 
$(\phi_{1}\otimes \phi_{1})\sigma (\jmath_{1} \omega \otimes u ) = 0$.
QED
\vspace{0.5cm}

Let $p : K^1 \otimes K^1
\longrightarrow K^2 $ be the multiplication map 
$ P_{1,1} = \phi_{2}\circ(\jmath_{1} \otimes \jmath_{1} )$
 and $p_{12} = p \otimes 1$.
We extend $\nabla $ to $\nabla _{1} $ and $\nabla_{2}$, respectively, by
defining 
$\nabla_{1} \equiv (\phi_{2} \otimes \phi_{1} ) 
   \circ D_{1} \circ (\jmath_{1} \otimes \jmath_{1} ) $
and
$\nabla_{2} \equiv (\phi_{1} \otimes \phi_{1} \otimes \phi_{1} ) 
   \circ D_{2} \circ (\jmath_{1} \otimes \jmath_{1} ). $
\vspace{0.5cm}

{\sc Proposition 6:} Let $ f\in A , \,\,\,\omega, \omega' \in K^1 $, and
    $\tau_{12} = \tau \otimes 1 $. Then

(1) $ \nabla_{1} ( \omega \otimes \omega' ) = 
   \delta \omega \otimes \omega' - p_{12}(\omega \otimes \nabla \omega' )$
    and

(2) $ \nabla_{2} (\omega \otimes \omega' ) = \nabla \omega \otimes \omega'
       + \tau_{12} (\omega \otimes \nabla \omega'),$
       if  $\sigma$ preserves $\triangle $.  
\vspace{0.5cm}

{\em Proof:} (1)  First, we observe that $\phi_{2} =
 p (\phi_{1}   \otimes \phi_{1} )$ from Lemma 5 and $\pi = 1 $
for the universal differential algebra $\Omega (A)$.
\begin{eqnarray*} 
\nabla_{1} (\omega \otimes \omega' ) 
&= & 
 (\phi_{2} \otimes \phi_{1} )(d\jmath_{1} \omega \otimes \jmath_{1}\omega') - 
        (\phi_{2} \otimes \phi_{1})
          (\jmath_{1}\omega \otimes D\jmath_{1} \omega' ) \\
& = & \delta \omega \otimes \omega' - p_{12}(\phi_{1} \otimes \phi_{1}
         \otimes  \phi_{1} )(\jmath_{1}\omega \otimes D\jmath_{1} \omega' )\\
& =& \delta \omega \otimes \omega' - p_{12}(\omega \otimes \nabla \omega' ).
\end{eqnarray*}
\begin{eqnarray*}
(2) \nabla_{2}(\omega \otimes \omega')
 &=& (\phi_{1} \otimes \phi_{1} \otimes \phi_{1} ) D_{2} 
       (\jmath_{1} \otimes \jmath_{1} )
          (\omega \otimes \omega' ) \\
& = & (\phi_{1} \otimes \phi_{1} \otimes \phi_{1} )
     (D\jmath_{1}\omega 
 \otimes \jmath_{1}\omega' + \sigma_{12} (\jmath_{1}\omega \otimes
                  D\jmath_{1}\omega')) \\
& = &\nabla \omega \otimes \omega'
+ (\phi_{1} \otimes \phi_{1} \otimes \phi_{1} ) 
   \sigma_{12}(\jmath_{1}\omega \otimes D\jmath_{1}\omega' ) \\
&= & \nabla \omega \otimes \omega'
+ (\phi_{1} \otimes \phi_{1} \otimes \phi_{1} ) 
   \sigma_{12}(\jmath_{1}\omega \otimes 
(\jmath_{1} \otimes \jmath_{1})\nabla\omega'
 + \jmath_{1} \omega \otimes u  ) 
\end{eqnarray*}
for some $u \in \Omega^1 \otimes_{A} \Omega^1 $
 such that $(\phi_{1} \otimes  \phi_{1})(u) = 0 $
since $(\phi_{1} \otimes \phi_{1} )( D\jmath_{1} \omega' 
 - (\jmath_{1} \otimes \jmath_{1} ) \nabla \omega' ) = 0$.
Now let us show that
$u \in \triangle $.
Let $u = \sum_{i}x_{i} \otimes y_{i}$ for $x_{i}, y_{i} \in \Omega^1 $.
Now  we can write $x_{i}$  uniquely as
$x_{i} = x_{i}^{(1)} + \jmath_{1}x_{i}^{(2)} $ 
    for $x_{i}^{(1)} \in ker \phi_{1}$
and $x_{i}^{(2)} \in K^1 $. Similarly for $y_{i}$.
Then $u = \sum_{i}\jmath_{1}x_{i}^{(2)} \otimes
                  \jmath_{1}y_{i}^{(2)} + \alpha $ 
for some $\alpha \in \triangle $.
From $(\phi_{1} \otimes \phi_{1} )(u) = 0$,
$\sum_{i}x_{i}^{(2)} \otimes y_{i}^{(2)} = 0$.
Hence $u = \alpha \in \triangle $.
Now that $\sigma (\triangle) \subset \triangle $,
$\sigma_{12}(\Omega^1 \otimes \triangle ) 
\subset \Omega^1 \otimes _{A} \Omega^1 \otimes_{A} ker \phi_{1} 
     + \triangle \otimes_{A} \Omega^1 $.
Thus 
$(\phi_{1} \otimes \phi_{1} \otimes \phi_{1})\sigma_{12}
(\jmath_{1} \omega \otimes u ) = 0$. QED
\vspace{0.5cm}

If we define the torsion  of $\nabla $ by $ T_{\nabla}
 = \delta - p \circ \nabla $, it is easy to see that  the necessary 
and sufficient condition for $T_{\nabla }$ to be $A$-bilinear
is that $p( 1 + \tau ) = 0 $ in the case where $\sigma $
preserves $\triangle$.
From now on, let us keep the $A$-bilinearity of $T$
as in the previous subsection. Thus $\sigma = -1 $.
Then $\tau = -1 $ and $T_{\nabla }$ is also $A$-bilinear. In this case,
we have $\nabla_{1} = p_{12} \nabla_{2} + T_{\nabla} \otimes 1 $.
If  $T = 0$, 
 $T_{\nabla} = \phi_{2} \circ d \circ \jmath_{1} 
       - p \circ (\phi_{1} \otimes \phi_{1} )
  \circ D \circ \jmath_{1} = (\phi_{2} - p\circ
    (\phi_{1} \otimes \phi_{1} ))\circ d \circ \jmath_{1} = 0 $.  
Now let us calculate explicitly connections $\nabla $ and 
curvatures $\nabla_{1} \circ \nabla $, $\nabla_{2} \circ \nabla $ on graphs
for the general nonzero $A$-bilinear torsion case.

As $D_{1}\circ D$ for complete symmetric graphs,  $\nabla_{1} \circ \nabla $  
is left $A$-linear and  not right $A$-linear.   
We note that
 $D\jmath_{1} e_{ij} = de_{i} \otimes de_{j} + 
\sum_{a} (1 + \Gamma_{ij}^{iaj} )e_{iaj}$ 
for $e_{ij} \in K^1 $ from Lemma 6.
Then it follows that 
\begin{eqnarray}
\nabla e_{ij} &=& (\phi_{1} \otimes \phi_{1}) D  \jmath_{1} e_{ij}
                    \nonumber \\
            &=& (\phi_{1} \otimes \phi_{1})(de_{i} \otimes de_{j} 
                       + \sum_{a }
                     (1 + \Gamma_{ij}^{iaj} ) e_{iaj} ) \\
     &=& \{\sum_{k}\phi_{1}(e_{ki} -e_{ik}) 
             + \sum_{a \neq j}(1 + \Gamma_{ij}^{iaj} ) 
              \phi_{1}(e_{ia})\} \otimes 
          \sum_{l}\phi_{1}(e_{lj} - e_{jl} ). \nonumber  
\label{eq:334}
\end{eqnarray}

For simplicity, we shall consider examples for the case where
$ 1 +  \Gamma_{ij}^{iaj} = \gamma $ for all $i, j \,\,\, ( i \neq j )$
 and all $a (\neq i, j) $. Thus 
\begin{eqnarray}
\nabla e_{ij} 
     = \{\sum_{k}\phi_{1}(e_{ki} -e_{ik}) 
             + \gamma \sum_{a \neq j} 
              \phi_{1}(e_{ia})\} \otimes 
          \sum_{l}\phi_{1}(e_{lj} - e_{jl} ).   
\label{eq:335}
\end{eqnarray}
\vspace{0.5cm}

{\sc Example 4:} Let ${\cal V} $ be a set $\{ x_{1}, x_{2} \} $ of two points
 as in the Connes-Lott's model$^{3}$. If $K^1 = \Omega^1 $, 
$D_{1}D = D_{2}D = 0$ from Corollary 3.  
If  $K^1 $ is the space 
generated by $\{ e_{12} \} $, $K^1 \otimes _{A} K^1 = 0 $ since
 $ e_{12} \otimes e_{12} = 0$.  Thus
 $\nabla_{1} \nabla = \nabla_{2} \nabla = 0 $.
The graph of two points seems to be too simple to have a nonzero curvature.
\vspace{0.5cm}

{\sc Example 5:} Let ${\cal V} $ be a set $\{x_{1}, x_{2}, x_{3} \}$
of three points and let $K^1 $ be the space generated by
$\{ e_{12}, e_{21}, e_{23}, e_{32} \}. $
Then $ ker \phi_{1} = < \{e_{13}, e_{31} \}>$,
$ker \phi_{2} = <\{e_{123}, e_{131}, e_{132},e_{213},e_{231},e_{312},
e_{313},e_{321} \}>$ and
$K^2 = <\{ e_{121}, e_{212}, e_{232}, e_{323} \}> $ with the obvious
splitting map. 
We note that
$\delta e_{12} = \phi_{2}d\jmath_{1} e_{12} = e_{121} + e_{212} 
  =  \delta e_{21}$, and
$\delta e_{23} = e_{232} + e_{323} =  \delta e_{32} $.  Moreover, 
\begin{eqnarray*}
\nabla e_{12} &=& (e_{21} - e_{12}) \otimes (e_{12} + e_{32}
                     -e_{21} - e_{23}) = e_{121} + e_{123} + e_{212} \\
\nabla e_{21} &=& (e_{12} + e_{32} - e_{21} - e_{23} + \gamma e_{23})
                    \otimes (e_{21} - e_{12})
                     = e_{121} + e_{212} + e_{321} 
\end{eqnarray*}
We can also calculate $\nabla e_{23}, \nabla e_{32} $ in the same way.
$ \nabla e_{23} = e_{123} + e_{232} + e_{323} $ and
 $ \nabla e_{32} = e_{232} + e_{321} + e_{323} $.
 
Now we have 
\begin{eqnarray*}
\nabla_{1}  \nabla e_{12}
    &=& \nabla_{1} ( e_{12} \otimes e_{21} + e_{21} \otimes e_{12} 
                  + e_{12} \otimes e_{23} ) \\
    &=& \delta e_{12} \otimes e_{21} -p_{12}(e_{12} \otimes \nabla e_{21})
       + \delta e_{21} \otimes e_{12} -p_{12}(e_{21} \otimes \nabla e_{12}) \\
        &&\mbox{} + \delta e_{12} \otimes e_{23} 
          -p_{12}(e_{12} \otimes \nabla e_{23}) \\
   &=& (e_{2121} - e_{1212}) + (e_{1212} - e_{2121} - e_{2123})
         +  e_{2123}  \\ 
   &=& 0,
\end{eqnarray*}
and
\begin{eqnarray*}
\nabla_{1}  \nabla e_{21}
    &=& \nabla_{1} ( e_{12} \otimes e_{21} + e_{21} \otimes e_{12} 
                  + e_{32} \otimes e_{21} ) \\
    &=& \delta e_{12} \otimes e_{21} -p_{12}(e_{12} \otimes \nabla e_{21})
       + \delta e_{21} \otimes e_{12} -p_{12}(e_{21} \otimes \nabla e_{12}) \\
      &&\mbox{}  +  \delta e_{32} \otimes e_{21} 
             - p_{12}(e_{32} \otimes \nabla e_{21}) \\
   &=& (e_{2121} - e_{1212}) + (e_{1212} - e_{2121} - e_{2123})
         +  e_{2321}  \\
   &=& e_{2321} - e_{2123}.
\end{eqnarray*}
On the other hand, $\nabla_{1} \nabla e_{23} = e_{2123} - e_{2321} $
and $\nabla_{1} \nabla e_{32} = 0 $.
If we factor out the second terms from  
$\nabla_{1} \nabla e_{21} $ and
$\nabla_{1} \nabla e_{23} $, we obtain an $A$-bilinear  curvature
$ Curv_{1}(e_{21}) = e_{2321}$ and $Curv_{1}(e_{23}) = e_{2123} $
 with
 $Curv_{1}(e_{12}) = Curv_{1}(e_{32}) = 0$.
 
Now we make use of  the following equation 
\[\nabla_{2}(e_{ijk}) = \nabla e_{ij} \otimes 
                  e_{jk} - e_{ij} \otimes \nabla e_{jk}\]
to obtain the curvature 
 $\nabla_{2}\circ \nabla $: 
\begin{eqnarray*}
\nabla_{2}\nabla e_{12} & = & - e_{1232} + e_{3212} =  
     - \nabla_{2}\nabla e_{32}, \\
\nabla_{2}\nabla e_{21} & = & - e_{2123} + e_{2321} =
       - \nabla_{2}\nabla e_{23}. 
\end{eqnarray*}

If the unwanted terms are factored out
from $p_{12}\nabla_{2}\nabla e_{ij}$, we have $A$-bilinear curvatures:
\begin{eqnarray*}
Curv_{2}(e_{12}) & = & 0, \hspace{1cm} Curv_{2}(e_{21}) = e_{2321} \\
Curv_{2}(e_{32}) & = & 0, \hspace{1cm}  Curv_{2}(e_{23}) = e_{2123}.
\end{eqnarray*}

If we add $e_{31}$ to $K^1$, not only the symmetry between the vertices 1 and
3 is broken but also the curvature depends on the parameter $\gamma $ as seen
in the next example. 
\vspace{0.5cm}

{\sc Example 6:} Let ${\cal V} $ be a set $\{x_{1}, x_{2}, x_{3} \}$
of three points and let $K^1 $ be the space generated by
$\{ e_{12}, e_{21}, e_{23}, e_{31}, e_{32} \}. $
Then $ ker \phi_{2}  $ = $ <\{e_{123}, e_{131},
 e_{132}, e_{213}, e_{313} \}>$ and
$K^2 = <\{ e_{121}, e_{212}, e_{231}, e_{232}, e_{312},
e_{321}, e_{323} \}> $. We note that
\begin{eqnarray*}
\delta e_{12} & = & e_{121} + e_{212} + e_{312}, \\
\delta e_{21} & = & e_{121} + e_{321} - e_{231} + e_{212}, \\
\delta e_{23} & = & e_{323} + e_{231} + e_{232}, \\
\delta e_{31} & = & e_{231} - e_{321} + e_{312}, \\
\delta e_{32} & = & e_{232} - e_{312} + e_{321} + e_{323},  
\end{eqnarray*}
and
\begin{eqnarray*}
\nabla e_{12} &=& e_{212} + e_{312} + e_{121} + e_{123}, \\
\nabla e_{21} &=& e_{121} + e_{321} + e_{212} - e_{231} 
                    + \gamma e_{231}, \\
\nabla e_{23} &=& e_{123} + e_{323}  + e_{232} + e_{231}, \\
\nabla e_{31} &=& e_{231}  - e_{321} + e_{312}
                + \gamma e_{321}, \\
\nabla e_{32} &=& e_{232} + e_{321} + e_{323} - e_{312} 
                   + \gamma e_{312}. 
\end{eqnarray*}

Then it follows that
\begin{eqnarray*}
\nabla_{1}\nabla e_{12} &=& 0,\\
\nabla_{1}\nabla e_{21} &=& -e_{2123} - \gamma e_{2312} 
            + \gamma (2 - \gamma )e_{2321},\\
\nabla_{1}\nabla e_{23} &=& e_{2123} - \gamma (e_{2312} + e_{2321}),\\
\nabla_{1}\nabla e_{31} &=& -e_{3123}   -\gamma e_{3121}
              + \gamma (2-\gamma )e_{3231}  - \gamma e_{3212},\\
\nabla_{1}\nabla e_{32} &=& -\gamma (e_{3231} +  e_{3212} 
                                         + e_{3121} + e_{3123}).
\end{eqnarray*}

Thus by factoring out the unwanted terms, we obtain an $A$-bilinear curvature.
\begin{eqnarray*}
Curv_{1}(e_{12}) &=& 0, \\
Curv_{1}(e_{21}) &=&  \gamma (2-\gamma )e_{2321},\\
Curv_{1}(e_{23}) &=& e_{2123}, \\
Curv_{1}(e_{31}) &=&  -\gamma e_{3121}  + \gamma (2 - \gamma )e_{3231},\\
Curv_{1}(e_{32}) &=& -\gamma e_{3212}.
\end{eqnarray*}

Similarly, we have
\begin{eqnarray*}
\nabla_{2}\nabla e_{12} &=& -e_{1232} + \gamma (e_{2312} + 
                                e_{3212} - e_{1231}),\\
\nabla_{2}\nabla e_{21} &=& -e_{2123} + \gamma e_{3121} 
            + \gamma (2 - \gamma )e_{2321},\\
\nabla_{2}\nabla e_{23} &=& e_{2123} + \gamma (e_{3123} -e_{2312} - e_{2321}),\\
\nabla_{2}\nabla e_{31} &=& e_{1231} - e_{3123}  
                         + \gamma (\gamma - 2) e_{3121}
                        + \gamma (2-\gamma )e_{3231},\\
\nabla_{2}\nabla e_{32} &=& e_{1232} -\gamma e_{3231} -\gamma  e_{3212} 
\end{eqnarray*}
and
\begin{eqnarray*}
Curv_{2}(e_{12}) &=& 0, \\
Curv_{2}(e_{21}) &=&  \gamma (2-\gamma )e_{2321},\\
Curv_{2}(e_{23}) &=& e_{2123}, \\
Curv_{2}(e_{31}) &=&
             \gamma (\gamma - 2 )e_{3121} + \gamma (2 - \gamma )e_{3231},\\
Curv_{2}(e_{32}) &=& -\gamma e_{3212}.
\end{eqnarray*}
\vspace{0.5cm}

We have one-parameter family of
 connections even on a subgraph of a complete symmetric graph.  
It is worthy to notice that $Curv_{1} = Curv_{2} $ 
for $\gamma = 1 $ in the above examples 
as in  complete symmetric graphs.
We note that
even though the torsion $T$ is free on a complete symmetric graph
(and hence $T_{\nabla } = 0$),
curvatures  do not vanish in general  on its  subgraphs as expected. 
\vspace{1cm}

{\bf \large IV. CONCLUSIONS}
\vspace{0.5cm}

Quantum spaces arise as models for
the description of the small scale structure of spacetime.
In this work, we have treated  graphs  as quantum spaces
and constructed their differential algebras by extending the formulation
of Dimakis et al$^{6}$ in such a manner that the calculation
of connections and curvatures can be done.
We have shown explicitly  the form of linear connections and
 curvatures of a given   complete
symmetric graph for the general nonzero torsion case. 
Also   $A$-bilinear curvatures 
 have been obtained for graphs
by factoring out the unwanted terms from  curvatures 
which are not $A$-bilinear.
An interesting question for further study is whether or not 
$Curv_{1} = Curv_{2} $  whenever $\gamma = 1 $  for any graphs.

There is one-parameter family of connections.
This fact parallels those for  other models
such as quantum planes and matrix geometries.
There is a metric, 
but no metric-compatible connection in general except the
complete symmetric graph 
with two vertices.
 Even though the torsion vanishes  and hence the curvature is zero 
on a complete symmetric graph,
a subgraph  obtained from it
by deleting some arrows can have a nonzero curvature 
in general as expected.
\vspace{1cm}

{\bf \large ACKNOWLEDGMENTS}
\vspace{0.5cm}

We thank F. M\"{u}ller-Hoissen for some comments and pointing out to us
Ref.\ 14.  One of the authors (S. Cho) is grateful to J. Madore
for a helpful discussion on metric.   
This work was supported by
Ministry of Education, Project No.\ BSRI-95-2414 and a grant from
TGRC-KOSEF 1995.

\vspace{2cm}

 \hspace{-1cm}$^{1}$ A. Connes,  {\em Non-commutative geometry } 
          (Academic, New York, 1994).\\

\hspace{-1cm}$^{2}$ S. L. Woronowicz,
           "Compact matrix pseudogroups,"
            Comm. Math. Phys. {\bf 111}, 613 (1987).\\

\hspace{-1cm}$^{3}$  A. Connes  and J.  Lott, 
            "Particle physics and non-commutative geometry,"
              Nucl. Phys. B (Proc. Suppl.)  {\bf 18},  29 (1990).\\

\hspace{-1cm}$^{4}$ J. Ambj$\o$rn,
         "Quantization of geometry," Lectures given at Les Houches,
           Session LXII, Fluctuating Geometries in Statistical Mechanics
             and Field Theory (1994). \\

\hspace{-1cm}$^{5}$ F. David,
         "Simplicial quantum gravity and random lattices,"
           Lectures given at Les Houches,
           Session LVII, Gravitation and Quantizations (1992). \\

\hspace{-1cm}$^{6}$ A. Dimakis  and F. M\"{u}ller-Hoissen,
        "Discrete differential calculus: Graphs, topologies, and
                gauge theory,"
             J. Math. Phys. {\bf 35}, 6703 (1994).\\

\hspace{-1cm}$^{7}$ A. Dimakis  and F.  M\"{u}ller-Hoissen,
         "Differential calculus and gauge theory on finite sets,"
             J. Phys. A. {\bf 27}, 3159 (1994).\\

\hspace{-1cm}$^{8}$  H.C.  Baehr, A. Dimakis and F. M\"{u}ller-Hoissen, 
            "Differential calculi on commutative algebras,"
             J. Phys. A. {\bf 28}, 3197 (1995).\\

\hspace{-1cm}$^{9}$ A. H. Chamseddine, G. Felder  and J.  Fr\"{o}hlich,
                "Gravity in non-commutative geometry,"
                  Comm. Math. Phys. {\bf 155}, 205 (1993).\\

\hspace{-1.1cm}$^{10}$ A. Sitarz,
          "Gravity from noncommutative geometry,"
        Class. Quant. Grav. {\bf 11}, 2127 (1994).\\

\hspace{-1.1cm}$^{11}$ A. Cuntz  and D. Quillen,
            "Algebra extensions and  nonsingularity," 
              J. Amer. Math. Soc. {\bf 8}, 251 (1995).\\

 \hspace{-1.1cm}$^{12}$ M. Dubois-Violette  and P. Michor,   
             "Connections on central bimodules in noncommutative 
                 differential geometry,"
              J. Geo. Phys. {\bf 20}, 218 (1996) \\

\hspace{-1.1cm}$^{13}$ M. Dubois-Violette and P.  Michor, 
          "D\'{e}rivations et calcul 
              diff\'{e}rentiel non commutatif II," 
            C. R. Acad. Sci. Paris S\'{e}rie I {\bf 319 }, 927 (1994).\\

\hspace{-1.1cm}$^{14}$ K. Bresser, F. M\"{u}ller-Hoissen, A. Dimakis  and
      A. Sitarz,  "A Noncommutative geometry of finite groups," 
             J. Phys. A. {\bf 29}, 2705 (1996).  \\

\hspace{-1.1cm}$^{15}$ J. Mourad,
                "Linear connections in non-commutative geometry,"          
                     Class. Quant. Grav. {\bf 12}, 965 (1995).\\

\hspace{-1.1cm}$^{16}$ N. Bourbaki, {\em Algebra I } 
          (Addison-Wesley Publishing Company, 1973).\\

\hspace{-1.1cm}$^{17}$ T. Brzezi\'{n}ski  and  S. Majid,
              "Quantum group, Gauge theory on quantum spaces," 
                       Comm. Math. Phys. {\bf 157}, 591 (1993).\\ 

\hspace{-1.1cm}$^{18}$ A. Gibbons, {\em Algorithmic graph theory} 
                    (Cambridge Univsity Press, 1985).\\

\hspace{-1.1cm}$^{19}$  M. Dubois-Violette,  J.  Madore, T.  Masson and 
                         J. Mourad,
                       "On curvature in noncommutative geometry," 
             J. Math. Phys. {\bf 37}, 4089 (1996). \\

\hspace{-1.1cm}$^{20}$ J. Madore, T. Masson  and J. Mourad, 
               "Linear connections on matrix geometries,"
                  Class. Quant. Grav. {\bf 12}, 1429 (1995).\\

\hspace{-1.1cm}$^{21}$ J. Madore, "Linear connections on fuzzy manifolds,"
             Class. Quant. Grav. {\bf 13}, 2109 (1996).\\

\hspace{-1.1cm}$^{22}$ A. Kehagias, J. Madore, J. Mourad and G. Zoupanos,
               "Linear connections on extended space-time,"
                  J. Math. Phys. {\bf 36}, 5855 (1995).\\

\hspace{-1.1cm}$^{23}$ M. Dubois-Violette, J. Madore, T. Masson and J. Mourad, 
                   "Linear connections on the quantum plane,"
                Lett. Math. Phys. {\bf 35}, 351 (1995).\\


\end{document}